\documentclass[ECP]{ejpecp} % replace ECP by EJP if needed.

\usepackage[utf8]{inputenc}

\newcommand{\Cov}[0]{\text{Cov}}
\newcommand{\Var}[0]{\text{Var}}

\newcommand{\diag}[0]{\text{diag}}

\newcommand{\R}[0]{\mathbb{R}}

\newcommand{\Prob}[0]{\mathds{P}}

\theoremstyle{remark}
\newtheorem{thm}{Theorem}[section]
\newtheorem{lem}[thm]{Lemma}

\newtheorem{rem}[thm]{Remark}

\renewcommand{\leq}{\leqslant} 
\renewcommand{\geq}{\geqslant}

\newcommand{\rd}{\mathrm{d}}

\newcommand{\cL}{\mathcal{L}}

\newcommand{\vW}{\mathbf{W}}
\newcommand{\vX}{\mathbf{X}}

\newcommand{\mvTheta}{\boldsymbol{\Theta}}

\newcommand{\bR}{\mathbb{R}}
\newcommand{\bS}{\mathbb{S}}

\DeclareMathOperator{\E}{\mathds{E}}

\DeclareMathOperator{\argmax}{argmax}

\DeclareMathOperator{\tr}{tr} 
  
\SHORTTITLE{Interacting particle systems}

\TITLE{Maximum likelihood estimation of potential energy in interacting particle systems from single-trajectory data}
\AUTHORS{%
  Xiaohui~Chen\footnote{University of Illinois at Urbana-Champaign, United States of America.
    \EMAIL{xhchen@illinois.edu}}}
 % Separate items with ;

\AMSSUBJ{60H15; 62M05.} % Edit. Separate items with ;
\VOLUME{0}
\YEAR{2020}
\PAPERNUM{0}
\DOI{10.1214/YY-TN}

\ABSTRACT{This paper concerns the parameter estimation problem for the quadratic potential energy in interacting particle systems from continuous-time and single-trajectory data. Even though such dynamical systems are high-dimensional, we show that the vanilla maximum likelihood estimator (without regularization) is able to estimate the interaction potential parameter with optimal rate of convergence simultaneously in mean-field limit {\it and} in long-time dynamics. This to some extend avoids the curse-of-dimensionality for estimating large dynamical systems under symmetry of the particle interaction.}

\begin{document}

%%%%%%%%%%%%%%%%%%%%%%%%%%%%%%%%%%%%%%%%%%%%%%%%%%%%%%%%%%%%%%%%%%%
%%                                                               %%
%% No need for \maketitle.                                       %%
%%                                                               %%
%%%%%%%%%%%%%%%%%%%%%%%%%%%%%%%%%%%%%%%%%%%%%%%%%%%%%%%%%%%%%%%%%%%

%%%%%%%%%%%%%%%%%%%%%%%%%%%%%%%%%%%%%%%%%%%%%%%%%%%%%%%%%%%%%%%%%%%
%%                                                               %%
%% Please replace what follows by the body of your article       %%
%% (up to the bibliography):                                     %%
%%                                                               %%
%%%%%%%%%%%%%%%%%%%%%%%%%%%%%%%%%%%%%%%%%%%%%%%%%%%%%%%%%%%%%%%%%%%

\section{Introduction}
\label{sec:introduction}

Dynamical systems of interacting particles have a wide range of applications on modeling collective behaviors in physics~\cite{PhysRevLett.96.104302}, biology~\cite{Mogilner:1999aa,Topaz:2006aa}, social science~\cite{MotschTadmor2014}, and more recently in machine learning as useful tools to understand the stochastic gradient descent (SGD) dynamics on neural networks~\cite{MeiE7665}. Due to the large number of particles, such dynamical systems are high-dimensional even for single-trajectory data from each particle,
and statistical learning problems for lower-dimensional interaction functionals from data are usually challenging~\cite{Bongini2017,LiLuMaggioniTangZhang2019,Lu14424}. In this paper, we propose a likelihood based inference to estimate the parameters of the interaction function induced by a potential energy in an $N$ interacting particle system based on their observed (continuous-time) trajectories, and establish its statistical guarantees.

\subsection{Interacting $N$-particle systems}

In statistical mechanics, microscopic behaviors of a large number of random particles are related to explain macroscopic physical quantities (such as temperature distributions). Specifically, a system of $N$ interacting particles $(X_{t}^{N,1},\dots,X_{t}^{N,N})$ can be described by stochastic differential equations (SDEs) of the form:
\begin{equation}
\label{eqn:N-particle-system_SDE}
    \rd X_{t}^{N,i} = {1 \over N} \sum_{j=1}^{N} b(X_{t}^{N,j}-X_{t}^{N,i}) \rd t + \sigma \rd W_{t}^{i},
\end{equation}
where $b : \bR^{d} \to \bR^{d}$ is a vector field representing the pairwise interaction between the particles, $(W_{t}^{1})_{t \geq 0}, \dots, (W_{t}^{N})_{t \geq 0}$ are $N$ independent copies of the standard Brownian motion in $\bR^{d}$ such that $W^{i}_{0} = 0$, and $\sigma \in \bR_{+}$ is a diffusion parameter which is assumed to be constant and known. The scaling in~\eqref{eqn:N-particle-system_SDE} puts us in the mean-field regime, where the pairwise interaction effect is weak and decays on the order of $1/N$ as the number of particles $N \to \infty$. Thus the total interaction effect remains $O(1)$.

In this paper, we consider the estimation problem of the interaction function $b$ parametrized by a linear approximation $b(x) = \Theta x$ for some unknown $d \times d$ (symmetric) positive-definite matrix $\Theta \succ 0$:
\begin{equation}
\label{eqn:N-particle-system_SDE_linear}
    \rd X_{t}^{N,i} = \Theta (\overline{X}^{N}_{t} - X_{t}^{N,i}) \rd t + \sigma \rd W_{t}^{i},
    \end{equation}
and $\overline{X}^{N}_{t} = N^{-1} \sum_{j=1}^{N} X_{t}^{N,j}$. Interaction in stochastic system~\eqref{eqn:N-particle-system_SDE_linear} relates to the {\it Hookean} behavior for capturing the {\it linear elasticity} where the interaction force $b$ scales linearly with deformation distance (due to compression and stretch) in the direction from $X_{t}^{N,j}$ to $X_{t}^{N,i}$. This type of interaction is extensively used to study the large-scale and long-time dynamics of protein folding as an elastic mass-and-spring network of small C$\alpha$ atoms~\cite{HalilogluBaharErman1997,EyalYangBahar2006}.

It is a classical result~\cite{Sznitman1991} that the $N$-particle SDE system~\eqref{eqn:N-particle-system_SDE} in the mean-field limit admits a unique strong solution $(X_{t}^{N,i})_{t \geq 0}$ if $b$ is (globally) Lipschitz, which is the case for the stochastic system~\eqref{eqn:N-particle-system_SDE_linear}. Based on the observed continuous-time and single-trajectory data of the particle movement $(X_{t}^{N,1})_{t \geq 0},\dots,(X_{t}^{N,N})_{t \geq 0}$ in the interacting particle system~\eqref{eqn:N-particle-system_SDE_linear}, our main focus is to estimate the interaction parameter $\Theta$ in the potential energy.

Note that the first-order dynamical system~\eqref{eqn:N-particle-system_SDE_linear} evolves as stochastic gradient flows in $\bR^{d}$:
\begin{equation}
\label{eqn:N-particle-system_SDE_linear_SGF}
{\rd X_{t}^{N,i} \over \rd t} = {1 \over N} \sum_{j=1}^{N} \nabla V(X_{t}^{N,j}-X_{t}^{N,i}) + \sigma \xi^{i}_{t},
\end{equation}
which corresponds to a quadratic potential energy $V(x) = {1 \over 2} x^{T} \Theta x$ and $\xi^{i}_{t}$ are i.i.d. standard Gaussian random vectors in $\bR^{d}$. The left-hand side of~\eqref{eqn:N-particle-system_SDE_linear_SGF} is the observed velocity vector of particle $i$ at time $t$ and the right-hand side of~\eqref{eqn:N-particle-system_SDE_linear_SGF} is a linear function of all particle trajectories at time $t$ corrupted by independent additive Gaussian noise. Thus, there are $Nd$ SDEs with observed $N$ trajectory data in dimension $d$ to solve in~\eqref{eqn:N-particle-system_SDE_linear_SGF} and estimation problem for $\Theta$ can be recast as a high-dimensional linear regression problem in an augmented space $\bR^{Nd}$ (cf. the equivalent form~\eqref{eqn:N-particle-system_SDE_linear_vector_form} in Section~\ref{sec:likelihood-ratio}). Nevertheless, the trajectory data are {\it temporally dependent} samples  since the particles are interacting and dynamic, so that theoretical guarantees on estimating structured coefficients in high-dimensional linear models with i.i.d. samples are no long applicable in our context~\cite{MR2807761}. Moreover, due to the {\it symmetry} of the particles in law, the regression coefficients have very special replicated block diagonal structure in $\bR^{Nd}$, which suggests that regularization techniques may not necessarily needed in our problem. Indeed, we show that a direct likelihood-ratio method suffices to estimate $\Theta$ with the optimal rate of convergence in this work.

\subsection{Stochastic Vlasov equation: decoupled mean-field limit}
\label{sec:Vlasov_equation}

Let $\rho_{t}^{N} = N^{-1} \sum_{j=1}^{N} \delta_{X_{t}^{N,j}}$ be the empirical measure of the $N$ particles at time $t$. Then we can alternatively write~\eqref{eqn:N-particle-system_SDE} as
\begin{equation}
\label{eqn:N-particle-system_SDE_2}
    \rd X_{t}^{N,i} = \Big( \int_{\bR^{d}} b(y-X_{t}^{N,i}) \rho_{t}^{N}(\rd y) \Big) \rd t + \sigma \rd W_{t}^{i},
\end{equation}
where the drift coefficient vector depends on the individual state $X_{t}^{N,i}$ and the distribution $\rho_{t}^{N}$ (due to interaction). As $N \to \infty$ (i.e., in the mean-field limit), the interaction contributed by any pair of particles in the $N$-particle system~\eqref{eqn:N-particle-system_SDE} vanishes and all particles are (asymptotically) i.i.d. since they have the same drift and diffusion coefficients driven by independent standard Brownian motion in $\bR^{d}$. By the law of large numbers, we see that for any fixed $t$, $\rho_{t}^{N} \to \rho_{t}$ as $N \to \infty$, and the dynamic system~\eqref{eqn:N-particle-system_SDE_2} becomes
\begin{equation}
\label{eqn:McKean-Vlasov-eq_N_copy}
    \rd Y_{t}^{i} = \Big( \int_{\bR^{d}} b(y-Y_{t}^{i}) \rho_{t}(\rd y) \Big) \rd t + \sigma \rd W_{t}^{i},
\end{equation}
where $(\rho_{t})_{t \geq 0}$ is a non-random measure flow. Since the particles are symmetric in distribution, $\rho_{t}$ is actually the limiting law of each particle $X_{t}^{i}, i =1,\dots,N$. This defines a system of independent {\it stochastic Vlasov equations}~\eqref{eqn:McKean-Vlasov-eq_N_copy} with $\rho_{t} = \cL(Y_{t}^{i})$, which is a class of Markov processes with nonlinear dynamics~\cite{McKean1966}.
%\begin{equation}
%\label{eqn:McKean-Vlasov-eq}
%    \rd Y_{t}^{i} = \Big( \int_{\bR^{d}} b(y-Y_{t}^{i}) \rho_{t}(y) \Big) \rd t + \sigma \rd W_{t}^{i}, \quad \mbox{where } ,
%\end{equation}

For quadratic potential $V(x) = {1 \over 2} x^{T} \Theta x$, the stochastic system~\eqref{eqn:N-particle-system_SDE_linear} depends on the empirical distribution $\rho_{t}^{N}$ through the empirical mean $\overline{X}_{t}^{N}$ and it can be decoupled and approximated by a system of independent mean-reverting processes. Averaging the $N$ SDEs in~\eqref{eqn:N-particle-system_SDE_linear}, we see that once again by the law of number numbers,
%of more general nonlinear and interacting $N$ particles~\eqref{eqn:N-particle-system_SDE} with a Lipschitz drift vector field $\fb$. Suppose that $X^{N,i}_{0} = \xi^{i}$ for some mean-zero i.i.d. random vectors $\xi^{i}$ in $\bR^{d}$. Averaging the $N$-particle system~\eqref{eqn:N-particle-system_SDE_quadratic}, we see that by the law of number numbers,
\[
\overline{X}^{N}_{t} = \overline{W}^{N}_{t} := {1 \over N} \sum_{i=1}^{N} W^{t}_{t} \to 0, \quad \mbox{as } N \to \infty,
\]
which means that the interaction effect $\overline{X}^{N}_{t}$ becomes deterministic and it is nicely {\it decoupled} in the mean-field limit. Thus we expect that the $i$-th particle process $(X_{t}^{N,i})_{t \geq 0}$ can be (independently) approximated by a limiting process $(Y^{i}_{t})_{t \geq 0}$ given by
\begin{equation}
\label{eqn:limiting_proc_OU}
\rd Y^{i}_{t} = -\Theta Y^{i}_{t} \rd t + \sigma \rd W^{i}_{t}.
\end{equation}
Note that the processes $(Y^{i}_{t})_{t \geq 0}$~\eqref{eqn:limiting_proc_OU} are {\it independent} copies of the Ornstein-Uhlenbeck (OU) processes, which is a linear dynamic system of $N$ independent particles.

\subsection{Existing literature}

Learnability (i.e., identifiability) of interaction functions in interacting particle systems under the coercivity condition were studied in~\cite{Lu14424,LiLuMaggioniTangZhang2019}. In the noiseless setting, estimation of the interaction kernel, a scalar-valued function of pairwise distance between particles in the system, was first studied in~\cite{Bongini2017} for single-trajectory data in the mean-field limit, where the rate of convergence is no faster than $N^{-1/d}$. To alleviate the curse-of-dimensionality, sparsity-promoting techniques were considered for some structured high-dimensional dynamical systems~\cite{Brunton3932,SchaefferTranWard2018}. \cite{Lu14424} showed that a least-squares estimator achieves the optimal rate of convergence (in the number of observed trajectories for each particle) for estimating the interaction kernel based on multiple-trajectory data sampled from a deterministic system with random initialization. Estimation of the diffusion parameter for interacting particle systems from noisy trajectory data was studied in~\cite{HuangLiuLu2019}. Consistency of parameter estimation of the general McKean-Vlasov equation by the maximum likelihood estimation is studied in~\cite{WEN2016237}. To the best of our knowledge, there is no existing work, regularized or not, establishes the optimal rate of convergence for interaction parameter estimation simultaneously in the large $N$ (mean-field limit) and large $t$ (long-time dynamics) regime. This work fills this gap for the linear elasticity interacting particle systems.

\subsection{Notation}

For two generic vectors $a,b \in \bR^{d}$, we use $a \cdot b = \sum_{j=1}^{d} a_{j} b_{j}$ to denote the inner product of $a$ and $b$. We use $\|a\| = (a \cdot a)^{1/2}$ to denote its Euclidean norm. For a generic matrix $M$, we use $\|M\|$ to denote its spectral norm. Denote the set of $d \times d$ (symmetric) positive-definite matrices by $\bS^{d \times d}_{+}$. We use $C$ and $c$ to denote positive universal constants whose values may vary from place to place.

\section{Maximum likelihood estimation}
\label{sec:likelihood-ratio}

We estimate $\Theta \succ 0$ by a likelihood based method. In view of propagation of chaos~\cite{Sznitman1991}, we assume that the initializations $X^{N,i}_{0}$ are some mean-zero independent random vectors in $\bR^{d}$. Let $\vW_{t} = (W^{1}_{t}, \dots, W^{N}_{t}) \in \bR^{Nd}$ be the $N$ stacked standard $d$-dimensional Brownian motion and $\vX^{N}_{t} = (X^{N,1}_{t}, \dots, X^{N,N}_{t}) \in \bR^{Nd}$ be the stacked observation process. Then system~\eqref{eqn:N-particle-system_SDE_linear} can be rewritten as a higher-dimensional mean-reverting process in the augmented space $\bR^{Nd}$:
\begin{equation}
\label{eqn:N-particle-system_SDE_linear_vector_form}
\rd \vX^{N}_{t} = - \mvTheta H \vX^{N}_{t} \rd t + \Sigma \rd \vW_{t},
\end{equation}
where $\mvTheta = \diag(\Theta,\dots,\Theta)$ is an $(Nd) \times (Nd)$ block diagonal matrix, $\Sigma = \diag(\sigma^{2},\dots,\sigma^{2})$ is an $(Nd) \times (Nd)$ diagonal matrix, and
\[
H = {1 \over N} \left(
\begin{array}{cccc}
(N-1) I_{d} & -I_{d} & \dots & -I_{d} \\
 -I_{d} & (N-1) I_{d} & \dots & -I_{d} \\
 \vdots & \vdots & \ddots & \vdots \\
  -I_{d} & -I_{d} & \dots & (N-1) I_{d} \\
\end{array}
\right)
\]
is an $(Nd) \times (Nd)$ {\it interaction matrix}. Note that $H$ is a projection matrix $H^{2} = H$, which implies that the interaction effect is homogeneous.

%For simplicity, we consider diagonal matrix $\Theta = \diag(\theta)$ where $\theta = (\theta_{1},\dots,\theta_{d}) \in \bR_{+}^{d}$ is the parameter of interest.

Let $P$ be the law of the standard $(Nd)$-dimensional Brownian motion $(\vW_{t})_{t \geq 0}$ and $Q$ be the law of the augmented observation process $(\vX^{N}_{t})_{t \geq 0}$. By the multivariate Girsanov theorem for changing measures (cf. Theorem 1.12 in \cite{Kutoyants2004}), the likelihood ratio of $(\vX^{N}_{t})_{t \geq 0}$ in~\eqref{eqn:N-particle-system_SDE_linear} and $(\vW_{t})_{t \geq 0}$ is given by the Radon-Nikodym derivative ${\rd Q \over \rd P} (\vX_{0}^{N,t}, \rho_{0}^{t}) =: e^{\overline{\ell}^{N}_{t}(A)}$, where
\begin{equation}
    \label{eqn:likelihood-ratio}
   \overline{\ell}_{t}^{N}(A) = \sum_{i=1}^{N} \Big[ -{1 \over 2} \int_{0}^{t} \|A(\overline{X}^{N}_{s} - X^{N,i}_{s})\|^{2} \, \rd s + \int_{0}^{t} A (\overline{X}^{N}_{s} - X^{N,i}_{s}) \cdot \rd X^{N,i}_{s} \Big],
\end{equation}
$A \in \bS^{d \times d}_{+}$, $\vX^{N,t}_{0} = (\vX_{s}^{N})_{s \in [0,t]}$, and $\rho_{0}^{t} = (\rho_{s})_{s \in [0,t]}$. Then the maximum likelihood estimator (MLE) for $\Theta$ is defined as
\begin{equation}
\label{eqn:MLE}
    \hat{\Theta}_{t}^{N} = \argmax_{A \in \bS^{d \times d}_{+}} \overline{\ell}_{t}^{N}(A).
\end{equation}

\begin{rem}[Computing the MLE]
The MLE $\hat{\Theta}_{t}^{N}$ in~\eqref{eqn:MLE} is a constrained optimization problem on the smooth manifold $\bS_+^{d \times d}$. Since the objective function $\overline{\ell}_{t}^{N}(A)$ is quadratic in $A$, we may easily maximize $\overline{\ell}_{t}^{N}(A)$ over all possible $d \times d$ matrices in closed form that may not necessarily be a symmetric positive-definite matrix (cf. equation~\eqref{eqn:MLE_formula} below), and then project the unconstrained maximizer $\tilde{\Theta}_{t}^{N} = \argmax_{A \in \R^{d \times d}} \overline{\ell}_{t}^{N}(A)$ into $\bS_+^{d \times d}$ by a procedure called {\it positive-definitization} that has been used to estimate high-dimensional covariance matrix (cf. Section 2.2 in \cite{ChenXuWu2013}). Let $\breve{\Theta}_{t}^{N} = \sum_{j=1}^d \breve{\lambda}_j u_j u_j^T$ be its eigen-decomposition of the symmetrized matrix $\breve{\Theta}_{t}^{N} = (\tilde{\Theta}_{t}^{N} + (\tilde{\Theta}_{t}^{N})^T) / 2$.

\begin{lem}[Positive-definitization preserves the rate]
\label{lem:positive-definitization}
If $\|\tilde{\Theta}_{t}^{N}-\Theta\| \leq r_n$, then $\| \check{\Theta}_{t}^{N} - \Theta \| \leq 3 r_n$, where $\check{\Theta}_{t}^{N} = \sum_{j=1}^d \max(\breve{\lambda}_j, r_n) u_j u_j^T$.

\end{lem}
\begin{proof}[Proof of Lemma~\ref{lem:positive-definitization}]
First, triangle inequality yields $\|\breve{\Theta}_{t}^{N}-\Theta\| \leq \|\hat{\Theta}_{t}^{N}-\Theta\| \leq r_n$. Note that
\begin{align*}
\|\check{\Theta}_{t}^{N} - \Theta \| \leq \|\check{\Theta}_{t}^{N}-\breve{\Theta}_{t}^{N}\| + \|\breve{\Theta}_{t}^{N}-\Theta\| \leq & \Big| \sum_{j=1}^d \big[ (\max(\breve{\lambda}_j, r_n) - \breve{\lambda}_j ) \big] u_j u_j^T \Big| + r_n \\
= & \max_{j \in [d]} \big| \max(\breve{\lambda}_j, r_n) - \breve{\lambda}_j \big| + r_n.
\end{align*}
If $\breve{\lambda}_j \leq 0$, then $|r_n - \breve{\lambda}_j| \leq r_n + |\breve{\lambda}_j| \leq r_n + |\breve\lambda_j - \lambda_j| \leq r_n + \|\breve{\Theta}_{t}^{N}  - \Theta \| \leq 2 r_n$, where $\lambda_j$ is the $j$-th positive eigenvalue of $\Theta$. If $\breve{\lambda}_j > 0$, then $|\max(\breve{\lambda}_j, r_n) - \breve{\lambda}_j| \leq r_n$. Thus, we have  $\|\check{\Theta}_{t}^{N} - \Theta \| \leq 3 r_n$.
\end{proof}

Lemma~\ref{lem:positive-definitization} states that the estimator $\check{\Theta}_{t}^{N} = \sum_{j=1}^d \max(\breve{\lambda}_j, r_n) u_j u_j^T$ projects back $\tilde{\Theta}_{t}^{N}$ into $\bS^{d \times d}$ with the same rate of convergence as $\tilde{\Theta}_{t}^{N}$. In view of the equivalent theoretical guarantee, we shall in practice compute the unconstrained version of the MLE $\tilde{\Theta}_{t}^{N}$ as our working definition of the MLE $\hat{\Theta}_{t}^{N}$. In Section~\ref{sec:rate}, we derive the rate of convergence for $\tilde{\Theta}_{t}^{N}$ and we shall use $\hat{\Theta}_{t}^{N}$ to mean $\tilde{\Theta}_{t}^{N}$. 
\end{rem}

Next we explain the intuition why the MLE~\eqref{eqn:MLE} works. Note that we can write the log-likelihood as
\begin{equation}
\label{eqn:MLE_equiv}
\overline{\ell}_{t}^{N}(A) = N \int_{0}^{t} \tr \Big[ M_{s} \Big(-{1 \over 2} A A^{T} + A \Theta \Big) \Big] \, \rd s + \sigma \sum_{i=1}^{N} \int_{0}^{t} A (\overline{X}^{N}_{s} - X^{N,i}_{s}) \cdot \rd W^{i}_{s},
\end{equation}
where $M_{s} =  N^{-1} \sum_{i=1}^{N} (\overline{X}^{N}_{s} - X^{N,i}_{s}) (\overline{X}^{N}_{s} - X^{N,i}_{s})^{T}$ is the instantaneous mean-field covariance matrix at time point $s$.

As discussed earlier in Section~\ref{sec:introduction}, since the interaction among the $N$ particles in the mean-field regime is weak, we expect that those particles can be decoupled by their independent analogs. Let $(Y^{1}_{t})_{t \geq 0}, \dots, (Y^{N}_{t})_{t \geq 0}$ be {\it independent} copies of the OU processes defined in~\eqref{eqn:limiting_proc_OU}. Effectively we can view $(Y^{1}_{t})_{t \geq 0}, \dots, (Y^{N}_{t})_{t \geq 0}$ as a {\it decoupled system} of the $N$-particle system $(X^{N,1}_{t})_{t \geq 0},\dots,(X^{N,N}_{t})_{t \geq 0}$. Based on the decoupled processes, we can approximate $\overline{\ell}_{t}^{N}(A)$ by
\[
\overline{\tilde{\ell}}_{t}^{N}(A) := N \int_{0}^{t} \tr \Big[ \tilde{M}_{s} \Big(-{1 \over 2} A A^{T} + A \Theta \Big) \Big] \, \rd s - \sigma \sum_{i=1}^{N} \int_{0}^{t} A Y^{i}_{s} \cdot \rd W^{i}_{s},
\]
where $\tilde{M}_{s} =  N^{-1} \sum_{i=1}^{N} Y^{i}_{s} {Y^{i}_{s}}^{T}$. Decompose
\[
\overline{\ell}_{t}^{N}(A) = \underbrace{(\overline{\ell}_{t}^{N}(A) - \overline{\tilde{\ell}}_{t}^{N}(A))}_{\mbox{\scriptsize decoupling error}} + \underbrace{(\overline{\tilde{\ell}}_{t}^{N}(A) - \E[\overline{\tilde{\ell}}_{t}^{N}(A)])}_{\mbox{\scriptsize OU fluctuation error}} + \underbrace{\E[\overline{\tilde{\ell}}_{t}^{N}(A)]}_{\mbox{\scriptsize signal}}.
\]
Concentration bounds developed in Section~\ref{sec:tech_lemmas} allow us to control the decoupling and OU fluctuation errors around zero. Thus information useful for the estimation purpose comes from the signal part
\[
\E[\overline{\tilde{\ell}}_{t}^{N}(A)] = \tr \Big[  \Big( \int_{0}^{t}  \E[\tilde{M}_{s}] \rd s \Big) \Big(-{1 \over 2} A A^{T} + A \Theta \Big) \Big].
\]
%For instance, if the initial distribution $\xi^{i} = 0$, then Lemma~\ref{lem:OU_moments} in Appendix~\ref{sec:OU_proc} implies that
Since the matrix $\int_{0}^{t}  \E[\tilde{M}_{s}] \rd s$ is positive-definite, we see that the maximizer of $\E[\overline{\tilde{\ell}}_{t}^{N}(A)]$ is $A^{*} = \Theta$. This means that on the population level, the MLE equals to the true parameter. Combining this with the decoupling and OU fluctuation errors, we can obtain the rate of convergence for the MLE $\hat{\Theta}_{t}^{N}$ in~\eqref{eqn:MLE}.

\section{Rate of convergence}
\label{sec:rate}

In this section, we derive the rate of convergence for estimating $\Theta$ by the MLE $\hat{\Theta}^{N}_{t}$ in~\eqref{eqn:MLE} from the continuous-time and single-trajectory data for each particle. Below is the main result of this paper.

\begin{thm}
\label{thm:rate}
Let $(X^{N,1}_{t})_{t \geq 0}, \dots, (X^{N,N}_{t})_{t \geq 0}$ be the $d$-dimensional $N$-particle system defined in~\eqref{eqn:N-particle-system_SDE_linear} with i.i.d. initialization $X^{N,i}_{0} \sim N(0, D)$ for $i=1,\dots,N$. Let $\theta_1 \geq \dots \geq \theta_d > 0$ and $\tau_1 \geq \dots \geq \tau_d \geq 0$ be the ordered eigenvalues of $\Theta$ and $D$, respectively. If
\begin{equation}\
\label{eqn:rate_assumptions}
t \geq C \kappa \left( {1 \over \theta_d} + {\tau_1 \over \sigma^2} \right) \quad \mbox{and} \quad \kappa d \sqrt{\log(d/\varepsilon) \over N} \leq c
\end{equation}
for some universal constants $C$ and $c$ where $\kappa = \theta_1 / \theta_d$ is the condition number of $\Theta$, then we have with probability at least $1-14\varepsilon$,
\begin{equation}
\label{eqn:rate}
\big\| \hat\Theta^N_t - \Theta \big\| \leq C d \log\left({d \over \varepsilon}\right) \sqrt{ \Big( {\tau_1 \over \sigma^2 t} + 1\Big) {\kappa \theta_1 \over N t} }.
\end{equation}
\end{thm}

Theorem~\ref{thm:rate} is non-asymptotic and has several appealing features. 

First, the Gaussian initialization is not essential and it can be relaxed to any i.i.d. initialization with sub-Gaussian distributions with the $\psi_2$ norm controlled by $\tau_1$. The sub-Gaussian tail is necessary to obtain the exponential concentration rate in $N$, and the effect of initialization in the approximating OU processes decays exponentially fast in $t$ (cf.~\eqref{eqn:OU_proc_solution} in Section~\ref{sec:OU_proc}). Error bound in~\eqref{eqn:rate} also reflects the diminishing effect of initialization $\tau_1 /(\sigma^2 t) \to 0$ as $t \to \infty$.

Second, the conditions in~\eqref{eqn:rate_assumptions} are mild. We do not assume that the $N$-particle processes $(X^{N,1}_{t})_{t \geq 0}, \dots, (X^{N,N}_{t})_{t \geq 0}$ starts from the stationary distribution (which is a restrictive assumption), and the (continuous) time complexity of the particle trajectories is sharp in the following sense. Suppose $\Theta = \diag(\theta,\dots,\theta)$ is an isotropic interaction matrix (i.e., $\kappa = 1$) and $X^{N,1}_{0} = \dots = X^{N,N}_{0} = 0$ (i.e., $\tau_1 = 0$). Observe that the decoupled $N$ copies of the OU processes to approximate the dynamics of the interacting $N$-particle system have the equilibrium distribution as $N(0, \sigma^{2} (2\theta)^{-1} I_d)$, and in view of~\eqref{eqn:OU_proc_solution}, it takes at least $\Omega(\theta^{-1})$ time for $(X^{N,1}_{t})_{t \geq 0}, \dots, (X^{N,N}_{t})_{t \geq 0}$ mixing to the steady states (modulo small decoupling errors in $t$). In particular, if $\theta$ is closer to zero, then the log-likelihood ratio becomes flatter and thus larger $t$ is necessary to see the information from samples of the stationary distribution. On the other hand, if the processes start from the stationary distribution, then this trajectory time lower bound in $t$ is not needed to obtain~\eqref{eqn:rate}, provided that $d \sqrt{N^{-1} \log(d/\varepsilon)} \leq c$ which is the case whenever $\epsilon$ is not too small, e.g., $\varepsilon \geq d \exp(-c^2 N/d^2)$.
%We also comment that the lower bound $400$ on number of particles $N$ is not crucial -- it certainly can be improved. In addition, the closer initial distribution to the stationary distribution, the better improvement can be achieved.

Third, it is known that the MLE $\tilde\theta_t$ of a one-dimensional OU process $\rd Y_t = -\theta Y_t \rd t + \sigma \rd W_t$ for the parameter $\theta > 0$ has the exact rate of convergence in the sense that $\sqrt{t} (\tilde\theta_t - \theta)$ converges in distribution to $N(0,2\theta)$ as $t \to \infty$ (cf. Example 1.35 in~\cite{Kutoyants2004}). Contrast this with our result, we see the rate of convergence in~\eqref{eqn:rate} is rate-optimal in both $t$ and $N$, as well as $\theta_1$. Specifically, for fixed $d$, we can obtain from Theorem~\ref{thm:rate} the large $N$ (mean-field limit) and large $t$ (long-time dynamics) asymptotics as:
\begin{equation}
\label{eqn:asymptotic}
\|\sqrt{t}(\hat{\Theta}^{N}_{t} - \Theta)\| = O_{P} \left( \theta_1 / \sqrt{N} \right),
\end{equation}
provided that the $N$ particle processes start from a chaotic distribution, which agrees with the exact asymptotics of the one-dimensional OU process (or more generally, $d$-dimensional OU processes with isotropic $\Theta$). In view of the architecture of approximating the $N$-particle system by $N$ independent OU processes, the rate we derived in Theorem~\ref{thm:rate} is rate-optimal in $t$ and $N$, modulo small decoupling errors. We shall highlight that long-time dynamic behavior in~\eqref{eqn:asymptotic} as $t \to \infty$ cannot be obtained from the classical theory of the propagation of chaos (cf.~\cite{Sznitman1991}), where Gronwall's lemma (cf. Appendix 1 in~\cite{RevuzYor1991}) is typically used to control the decoupling error between the interacting $N$-particle processes and their independent analogs. In such case, the decoupling error is exponentially increasing in $t$, and thus it cannot be used to yield the rate $t^{-1/2}$. Our argument is tailored to the quadratic structure of the interacting potential $V(x) = {1 \over 2} x^{T} \Theta x$, which allows for a far more efficient decoupling strategy (cf. Lemma~\ref{lem:decoupling_error_bound}).

%Third, our concentration bound~\eqref{eqn:rate} allows the dimension of each particle process to increase at the speed $d = o(Nt)$ while maintaining a converging rate. All existing work in literature (\textcolor{red}{need citation here}) considers either $d = 1$ or $d \geq 1$ but fixed. In contrast, our estimator is statistically valid in terms of estimation in the increasing $d$ regime.

\begin{proof}[Proof of Theorem~\ref{thm:rate}.]
Since the objective function $\overline{\ell}_{t}^{N}(A)$ in~\eqref{eqn:MLE_equiv} is quadratic in $A$, the first-order optimality condition for the unconstrained optimization problem implies that the MLE of $\Theta$ satisfies
\begin{equation}
\label{eqn:MLE_formula}
\hat{\Theta}^{N}_{t} = \Theta + \Big(  \int_{0}^{t} M_{s} \rd s \Big)^{-1} \Big( \int_{0}^{t} {1 \over N} \sum_{i=1}^{N} \rd W_{s}^{i} \otimes (\overline{X}^{N}_{s} - X^{N,i}_{s}) \Big) \, \sigma,
%{{1 \over N} \sum_{i=1}^{N} \int_{0}^{t} (\overline{X}^{N}_{s,j} - X^{N,i}_{s,j}) \rd W^{i}_{s,j} \over {1 \over N} \sum_{i=1}^{N} \int_{0}^{t} | \overline{X}^{N}_{s,j} - X^{N,i}_{s,j} |^{2} \rd s}, \quad j = 1,\dots,d.
\end{equation}
where we recall $M_{s} = N^{-1} \sum_{i=1}^{N} (\overline{X}^{N}_{s} - X^{N,i}_{s}) (\overline{X}^{N}_{s} - X^{N,i}_{s})^{T}$ and $a \otimes b$ denotes the tensor product of two vectors $a$ and $b$, i.e., $(a \otimes b)_{jk} = a_j b_k$. Let $(Y_{t}^{i})_{t \geq 0}, i \in [N]$ be independent copies of the OU process driven by the same Brownian motion $(W^{i}_{t})_{t \geq 0}$ in $(X_{t}^{N,i})_{t \geq 0}$, namely,
\[
\rd Y_{t}^{i} = - \Theta Y_{t}^{i} \rd t + \sigma \rd W^{i}_{t} \quad \mbox{with} \quad Y^{i}_{0} = X^{N,i}_{0}.
\]
First, by triangle inequality,
\[
\Big\| {1 \over N} \sum_{i=1}^{N} \int_{0}^{t}  \rd W_{s}^{i} \otimes (\overline{X}^{N}_{s} - X^{N,i}_{s}) \Big\| \leq \Big\| {1 \over N} \sum_{i=1}^{N} \int_{0}^{t} \rd W_{s}^{i} \otimes (\overline{X}^{N}_{s} - X^{N,i}_{s} + Y^{i}_{s}) \Big\| + \Big\| {1 \over N} \sum_{i=1}^{N} \int_{0}^{t} \rd W_{s}^{i} \otimes Y^{i}_{s} \Big\| 
\]
Let $g(\varepsilon, N, d) = 1 + \sqrt{N^{-1} \log(d/\varepsilon)}$. Clearly, $g(\varepsilon, N, d) \leq g(\varepsilon, 1, d)$. By~\eqref{eqn:concentration_bound} in Lemma~\ref{lem:concentration_moment_bounds_OU} and~\eqref{lem:decoupling_error_bound_1} in Lemma~\ref{lem:decoupling_error_bound}, there exists a universal constant $C > 0$ such that for any $\varepsilon \in (0,1)$,
\begin{equation}
\label{eqn:MLE_formula_numerator}
\Big\| {1 \over N} \sum_{i=1}^{N} \int_{0}^{t}  \rd W_{s}^{i} \otimes (\overline{X}^{N}_{s} - X^{N,i}_{s}) \Big\| \leq C d g(\varepsilon, 1, 1) \sqrt{{(\tau_1+\sigma^2 t) \log(d/\varepsilon) \over N \theta_d}} \leq C d \log(d/\varepsilon) \sqrt{{\tau_1+\sigma^2 t \over N \theta_d}}
\end{equation}
holds with probability at least $1-8\varepsilon$. Next, since
\[
\Big\| \Big(  \int_{0}^{t} M_{s} \rd s \Big)^{-1} \Big\| = {1 \over \lambda_{\min}\Big(  \int_{0}^{t} M_{s} \rd s \Big)}= {1 \over \lambda_d\Big(  \int_{0}^{t} M_{s} \rd s \Big)},
\]
where $\lambda_j(M)$ is the $j$-th eigenvalue of a symmetric positive semidefinite matrix $M$ and $\lambda_{\min}(M)$ is the smallest eigenvalue of $M$, we shall derive a lower bound for  $\lambda_{\min}\Big(  \int_{0}^{t} M_{s} \rd s \Big)$. Note that
\begin{align*}
\int_{0}^{t} M_{s} \rd s = & \int_{0}^{t}  \E [Y^{i}_{s} {Y^{i}_{s}}^T] \rd s + {1 \over N} \sum_{i=1}^{N} \int_{0}^{t} ( Y^{i}_{s} {Y^{i}_{s}}^T - \E Y^{i}_{s} {Y^{i}_{s}}^T ) \rd s \\
& \qquad + {1 \over N} \sum_{i=1}^{N} \int_{0}^{t} [ (\overline{X}^{N}_{s} - X^{N,i}_{s}) (\overline{X}^{N}_{s} - X^{N,i}_{s})^T - Y^{i}_{s} {Y^{i}_{s}}^T ] \rd s.
\end{align*}
By Weyl's inequality, we have for all $j \in [d]$,
\begin{align*}
\Big| \lambda_j \Big(\int_{0}^{t} M_{s} \rd s\Big) - \lambda_j \Big(\int_{0}^{t}  \E [Y^{i}_{s} {Y^{i}_{s}}^T] \rd s\Big) \Big| \leq & \Big\| {1 \over N} \sum_{i=1}^{N} \int_{0}^{t} ( Y^{i}_{s} {Y^{i}_{s}}^T - \E Y^{i}_{s} {Y^{i}_{s}}^T ) \rd s \Big\| \\
& + \Big\| {1 \over N} \sum_{i=1}^{N} \int_{0}^{t} [ (\overline{X}^{N}_{s} - X^{N,i}_{s}) (\overline{X}^{N}_{s} - X^{N,i}_{s})^T - Y^{i}_{s} {Y^{i}_{s}}^T ] \rd s \Big\|.
\end{align*}
In addition, since $\E[Y^{i}_{s} {Y^{i}_{s}}^T] = \Cov(Y^i_s) = e^{-\Theta t} D e^{-\Theta t} + \sigma^2 (2\Theta)^{-1} (I_d - e^{-2\Theta t})$, by a second application of Weyl's inequality and Jensen's inequality, we have for all $j \in [d]$,
\begin{align*}
\Big| \lambda_j \Big(\int_0^t \E[Y^{i}_{s} {Y^{i}_{s}}^T] \rd s \Big) - \lambda_j \big(t \sigma^2 (2\Theta)^{-1} \big) \Big| \leq & \Big\| \int_0^t e^{-\Theta s} D e^{-\Theta s} \rd s - \int_0^t \sigma^2 (2\Theta)^{-1} e^{-2\Theta s} \rd s\Big\| \\
\leq & \int_0^t \big[ \tau_1 e^{-2s\theta_d} + \sigma^2 (2\theta_d)^{-1} e^{-2s\theta_d} \big] \rd s \leq {\tau_1 \over 2\theta_d} + {\sigma^2 \over 4\theta_d^2}.
\end{align*}
Combining the last two inequalities and using~\eqref{eqn:moment_bound} in Lemma~\ref{lem:concentration_moment_bounds_OU} and~\eqref{lem:decoupling_error_bound_2} in Lemma~\ref{lem:decoupling_error_bound}, we deduce that there is a universal constant $c$ such that with probability at least $1-6\varepsilon$,
\begin{equation}
\label{eqn:MLE_formula_denominator}
\lambda_{\min} \Big(\int_{0}^{t} M_{s} \rd s\Big) \geq {\sigma^2 t \over 2 \theta_1} - {\tau_1 \over 2\theta_d} - {\sigma^2 \over 4\theta_d^2} - C d {\tau_1 + \sigma^2 t \over 2 \theta_d} \sqrt{\log(d/\varepsilon) \over N} \geq c {\sigma^2 t \over \theta_1},
\end{equation}
where the last inequality follows from the conditions in~\eqref{eqn:rate_assumptions}. Now, putting together~\eqref{eqn:MLE_formula},~\eqref{eqn:MLE_formula_numerator} and~\eqref{eqn:MLE_formula_denominator}, we conclude that there exists a universal constant $C$ such that~\eqref{eqn:rate} holds with probability at least $1-14\varepsilon$.
\end{proof}

\section{Technical lemmas}
\label{sec:tech_lemmas}

This section provides key technical results for bounding the fluctuation of the OU process and the decoupling error of the $N$-particle system by the associated $N$ independent $d$-dimensional Ornstein-Uhlenbeck processes. 

\subsection{Concentration inequalities for Ornstein-Uhlenbeck processes}
\label{sec:OU_proc}

The $d$-dimensional Ornstein-Uhlenbeck (OU) process $(Y_{t})_{t \geq 0}$ with an isotropic diffusion parameter is a mean-reverting stochastic process defined by the following stochastic differential equation:
\begin{equation}
\label{eqn:OU_proc}
\rd Y_{t} = -\Theta Y_{t} \rd t + \sigma \rd W_{t},
\end{equation}
where $(W_{t})_{t \geq 0}$ is the standard Brownian motion in $\bR^d$ such that $W_{0} = 0$, $\Theta \in \bS^{d \times d}_+$ is the mean-reverting drift parameter, and $\sigma > 0$ is the scalar diffusion parameter. The solution to~\eqref{eqn:OU_proc} is given by
\begin{equation}
\label{eqn:OU_proc_solution}
Y_t = e^{-\Theta t} Y_0 + \sigma \int_0^t e^{\Theta (s-t)} \, \rd W_s,
\end{equation}
where $e^{\Theta} = \sum_{k=0}^\infty \Theta^k / k!$ is the matrix exponential function. From~\eqref{eqn:OU_proc_solution}, we see that $Y_t | Y_0 \sim N( e^{-\Theta t} Y_0, \sigma^2 \int_0^t e^{2\Theta(s-t)} \rd s)$. Suppose the initialization $Y_0$ has mean zero and covariance matrix $D$, and $Y_0$ is independent of $(W_{t})_{t \geq 0}$. Then we have $\E[Y_t] = 0$ and $\Sigma_t := \Cov(Y_t) = e^{-\Theta t} D e^{-\Theta t} + \sigma^2 (2\Theta)^{-1} (I_d - e^{-2\Theta t})$. As $t \to \infty$, the stationary distribution $Y_{t}$ is $N(0, \sigma^2 (2\Theta)^{-1})$.

Denote $g(\varepsilon, N, d) = 1 + \sqrt{N^{-1} \log(d/\varepsilon)}$. Let $\theta_1 \geq \dots \geq \theta_d > 0$ and $\tau_1 \geq \dots \geq \tau_d \geq 0$ be the ordered eigenvalues of $\Theta$ and $D$, respectively.

\begin{lem}[Concentration inequalities for the OU process]
\label{lem:concentration_moment_bounds_OU}
Let $\rd Y_{t} = -\Theta Y_{t} \rd t + \sigma \rd W_{t}$ be the $d$-dimensional OU process defined in~\eqref{eqn:OU_proc} with $Y_{0} \sim N(0,D)$ being independent of $(W_{t})_{t \geq 0}$. Suppose that $(Y^{1}_{t})_{t \geq 0}, \dots, (Y^{N}_{t})_{t \geq 0}$ are independent copies of $(Y_{t})_{t \geq 0}$. Then there exists a universal constant $C > 0$ such that for any $\varepsilon \in (0,1)$,
\begin{align}
\label{eqn:moment_bound}
& \Prob \Big( \Big\| {1 \over N} \sum_{i=1}^{N} \int_{0}^{t} ( Y^{i}_{s} {Y^{i}_{s}}^T - \E Y^{i}_{s} {Y^{i}_{s}}^T ) \rd s \Big\| \geq C d {\tau_1+\sigma^2 t \over \theta_d} \Big( {\log(d/\varepsilon) \over N} + \sqrt{\log(d/\varepsilon) \over N} \Big) \Big) \leq 2\varepsilon, \\
\label{eqn:concentration_bound}
& \Prob \Big( \Big\| {1 \over N} \sum_{i=1}^{N} \int_{0}^{t} \rd W^{i}_{s} \otimes Y^{i}_{s} \Big\| \geq C d g(\varepsilon, N, 1) \sqrt{{(\tau_1+\sigma^2 t) \log(d/\varepsilon) \over N \theta_d}} \Big) \leq 4 \varepsilon.
\end{align}
\end{lem}

\begin{proof}[Proof of Lemma~\ref{lem:concentration_moment_bounds_OU}]
Denote $U^{i}_{s} = Y^{i}_{s} {Y^{i}_{s}}^T - \E Y^{i}_{s} {Y^{i}_{s}}^T$ and $U^i_s = (U^i_{s,jk})_{j,k \in [d]}$. Let $\lambda > 0$ and fix a $j, k \in [d]$. By the exponential Markov inequality and the independence of the processes $(Y_t^{i})_{t \geq 0}$, we have for all $x > 0$,
\begin{equation}
\label{eqn:exponential_markov_ineq}
\Prob \Big( {1 \over N} \sum_{i=1}^N \int_{0}^{t} U^{i}_{s,jk} \rd s \geq x \Big) \leq \exp(-\lambda x) \prod_{i=1}^N \E\Big[\exp\Big( {\lambda \over N} \int_{0}^{t} U^{i}_{s,jk} \rd s \Big)\Big].
\end{equation}
%By Jensen's inequality and
%\[
%\E\Big[\exp\Big(\lambda \int_{0}^{t} {1 \over N} \sum_{i=1}^N U^{i}_{s,jk} \rd s \Big)\Big] \leq {1 \over t} \int_{0}^{t} \E \Big[\exp\Big(t \lambda {1 \over N} \sum_{i=1}^N U^{i}_{s,jk} \Big) \Big] \rd s =  {1 \over t} \int_{0}^{t} \prod_{i=1}^N \E \Big[\exp\Big({t \lambda \over N} U^{i}_{s,jk} \Big) \Big] \rd s.
%\]
Since $Y_0^i \sim N(0,D)$, we have $Y_t \sim N(0, e^{-\Theta t} D e^{-\Theta t} + \sigma^2 (2\Theta)^{-1} (I_d - e^{-2\Theta t}))$ for $t \geq 0$. Thus, for all $0 \leq s \leq t$,
\[
\Var(Y^i_{sj}) \leq \|\Cov(Y_s^i)\| \leq \|e^{-\Theta s}\|^2 \, \|D\| + \sigma^2 \, \|(2\Theta)^{-1}\| \, \|I_d - e^{-2\Theta s}\| \leq \tau_1 e^{-2s\theta_d} + \sigma^2 (2 \theta_d)^{-1},
\]
and $\|Y^i_{sj}\|_{\psi_2}^2 \leq C \, \Var(Y^i_{sj}) \leq C\, (\tau_1 e^{-2s\theta_d} + \sigma^2 (2 \theta_d)^{-1})$ for some universal constant $C$. By Lemma 2.7.7 in~\cite{Vershynin2018_Cambridge}, $\|Y_{sj}^i Y_{sk}^i\|_{\psi_1} \leq \|Y_{sj}^i\|_{\psi_2} \|Y_{sk}^i\|_{\psi_2} \leq C \, (\tau_1 e^{-2s\theta_d}+ \sigma^2 (2 \theta_d)^{-1})$. This together with Jensen's inequality imply that 
\[
\Big\| \int_0^t U^i_{s,jk} \rd s \Big\|_{\psi_1} \leq \int_0^t \left\| U^i_{s,jk} \right\|_{\psi_1} \rd s \leq C \, \theta_d^{-1} (\tau_1 + \sigma^2 t).
\]
Using the moment-generating function property of sub-exponential random variables (cf. Proposition 2.7.1 in~\cite{Vershynin2018_Cambridge}), we have
\[
\Big|{\lambda \over N} \Big| \leq {c \over \max_{i} \|\int_0^t U^i_{s,jk} \rd s\|_{\psi_1}}  \implies  \E \Big[\exp\Big({\lambda \over N} \int_0^t U^{i}_{s,jk} \rd s \Big) \Big] \leq \exp \Big( C {\lambda^2 \over N^2} \Big\|\int_0^t U^i_{s,jk} \rd s \Big\|_{\psi_1}^2 \Big).
\]
Thus, for $0 < \lambda < c N \theta_d / (\tau_1 + \sigma^2 t)$, we can bound
\[
 \E \Big[\exp\Big({\lambda \over N} \int_0^t U^{i}_{s,jk} \rd s \Big) \Big]  \leq  \exp \Big( C {\lambda^2 \over N^2} \sum_{i=1}^N \Big\|\int_0^t U^i_{s,jk} \rd s \Big\|_{\psi_1}^2 \Big) \leq \exp \Big( C {\lambda^2 (\tau_1 + \sigma^2 t)^2 \over N \theta_d^2} \Big).
\]
Combining the last inequality with~\eqref{eqn:exponential_markov_ineq} and optimizing the bound over $\lambda$, we have
\[
\Prob \Big( {1 \over N} \sum_{i=1}^N \int_{0}^{t} U^{i}_{s,jk} \rd s \geq x \Big) \leq \exp \Big\{ - \min \Big[ {N \theta_d^2 x^2 \over 4C (\tau_1+\sigma^2 t)^2}, \; {c N \theta_d x \over 2 (\tau_1+\sigma^2 t)}\Big] \Big\}.
\]
Equivalently, we can write the above inequality as: for all $\varepsilon \in (0,1)$,
\begin{equation}
\label{eqn:moment_bound_1d}
\Prob \Big( {1 \over N} \sum_{i=1}^N \int_{0}^{t} U^{i}_{s,jk} \rd s \geq C {\tau_1+\sigma^2 t \over \theta_d} \Big( {\log(1/\varepsilon) \over N} + \sqrt{\log(1/\varepsilon) \over N} \Big) \Big) \leq \varepsilon.
\end{equation}
Applying the same argument to $-U^{i}_{s,jk}$, together with the union bound over $j,k\in[d]$ and using the fact that $\|M\| \leq d \max_{j,k \in [d]} |M_{jk}|$ for any $d \times d$ matrix $M$, we obtain~\eqref{eqn:moment_bound}.

Next we prove~\eqref{eqn:concentration_bound}. Denote $Z^{i}_{t} = \int_{0}^{t} \rd W^{i}_{s} \otimes Y^{i}_{s}$. Clearly $\E[Z^{i}_{t}] = 0_{d \times d}$. Moreover, $(Z^{i}_{t})_{t \geq 0}, \dots, (Z^{N}_{t})_{t \geq 0}$ are independent continuous local martingales vanishing at zero with quadratic variation $[Z^{i}]_{t} := ([Z^{i}_{jk}]_{t})_{j,k\in[d]}$ is given by $[Z^{i}_{jk}]_{t} = \int_{0}^{t} |Y^{i}_{sk}|^{2} \rd s$. Moreover, the process $\overline{Z}^{N}_{s} = {1 \over N} \sum_{i=1}^{N} Z^{i}_{s}$ is also a local martingale vanishing at zero with quadratic variation $[\overline{Z}^{N}_{jk}]_{t} = {1 \over N^{2}} \sum_{i=1}^{N} \int_0^t |Y^{i}_{sk}|^2 \rd s$. Note that $\E[Z^{i}_{jk}]_{t} \leq \int_0^t \|\Sigma_s\| \rd s \leq (2\theta_d)^{-1} (\tau_1 + \sigma^2 t)$ and $\E [\overline{Z}^{N}_{jk}]_{t} \leq (2N\theta_d)^{-1} (\tau_1 + \sigma^2 t)$. Applying~\eqref{eqn:moment_bound_1d} on the diagonal entries $(U^i_{s,kk})_{k\in[d]}$, 
%we see that
%\begin{equation}
%\label{eqn:qv_concentraiton}
%\Prob \Big( [\overline{Z}^{N}_{jk}]_{t} - \E [\overline{Z}^{N}_{jk}]_{t} \geq  C {t \over N} (\tau_1 + {\sigma^2 \over 2 \theta_d}) \Big( {\log(1/\varepsilon) \over N} + \sqrt{\log(1/\varepsilon) \over N} \Big) \Big) \leq \varepsilon.
%\end{equation}
we have $\Prob([\overline{Z}^{N}_{jk}]_{t} \geq C t g(\varepsilon, N, 1)^2 (\tau_1 + \sigma^2 t) (N\theta_d)^{-1}) \leq \varepsilon$. By Bernstein's inequality for continuous local martingales (cf. Exercise (3.16) on page 145 in~\cite{RevuzYor1991}), we get
\begin{align*}
& \Prob\Big( \sup_{0 \leq s \leq t} {1 \over N} \sum_{i=1}^{N} Z^{i}_{s,jk} \geq C x g(\varepsilon, N, 1) \sqrt{\tau_1+\sigma^2 t \over N \theta_d} \Big) \\
\leq & \Prob\Big( \sup_{0 \leq s \leq t} \overline{Z}^{N}_{s,jk} \geq C x g(\varepsilon, N, 1) \sqrt{\tau_1+\sigma^2 t \over N \theta_d}, \, [\overline{Z}^{N}]_{t} \leq C g(\varepsilon, N, 1)^2 {\tau_1+\sigma^2 t \over N \theta_d} \Big) + \varepsilon \leq \exp \Big(- {C x^2 \over 2} \Big) + \varepsilon.
\end{align*}
Choosing $x = \sqrt{2\log(d^2/\varepsilon)}$ and applying the union bound, we have with probability at least $1 - 2 \varepsilon$,
\[
\max_{j,k\in[d]} \sup_{0 \leq s \leq t} {1 \over N} \sum_{i=1}^{N} Z^{i}_{s,jk} \leq C g(\varepsilon, N, 1) \sqrt{{ (\tau_1+\sigma^2 t) \log(d/\varepsilon) \over N \theta_d}}
\]
Applying the same argument to $-Z^{i}_{s}$ and using $\|M\| \leq d \max_{j,k \in [d]} |M_{jk}|$ for any $d \times d$ matrix $M$, we obtain~\eqref{eqn:concentration_bound}. 
\end{proof}

\subsection{Decoupling error bounds for $N$-particle systems}
\label{sec:decoupling_error_bound}

\begin{lem}[Decoupling error bounds for the $N$-particle system]
\label{lem:decoupling_error_bound}
Let $(X^{N,1}_{t})_{t \geq 0}, \dots, (X^{N,N}_{t})_{t \geq 0}$ be the $d$-dimensional $N$-particle system defined in~\eqref{eqn:N-particle-system_SDE_linear} with i.i.d. initialization $X^{N,i}_{0} \sim N(0, D)$ for $i=1,\dots,N$. Then there exists a universal constant $C > 0$ such that for any $\varepsilon \in (0,1)$,
\begin{align}
\label{lem:decoupling_error_bound_1}
& \Prob \Big( \Big\| {1 \over N} \sum_{i=1}^{N} \int_{0}^{t} \rd W^{i}_{s} \otimes (\overline{X}^{N}_{s} - X^{N,i}_{s}+Y^{i}_{s}) \Big\| \geq C d g(\varepsilon, 1, 1) \sqrt{ (\tau_1 + \sigma^2 t) \log(d/\varepsilon) \over N \theta_d} \Big) \leq 4 \varepsilon, \\ \notag
%& \Prob \Big( \Big\|{1 \over N} \sum_{i=1}^{N} \int_{0}^{t} [ (\overline{X}^{N}_{s} - X^{N,i}_{s}) (\overline{X}^{N}_{s} - X^{N,i}_{s})^T - Y^{i}_{s} {Y^{i}_{s}}^T ] \rd s \Big\| \leq C d t (\tau_1 + {\sigma^2 \over 2 \theta_d}) \sqrt{\log(d/\varepsilon) \over N} \Big) \leq 4\varepsilon.
& \Prob \Big( \Big\|{1 \over N} \sum_{i=1}^{N} \int_{0}^{t} [ (\overline{X}^{N}_{s} - X^{N,i}_{s}) (\overline{X}^{N}_{s} - X^{N,i}_{s})^T - Y^{i}_{s} {Y^{i}_{s}}^T ] \rd s \Big\| \\
\label{lem:decoupling_error_bound_2}
& \qquad \qquad \qquad \qquad \qquad \qquad \qquad \leq C d {g(\varepsilon,1,d) \over \sqrt{N}} \Big[ {g(\varepsilon,1,d) \over \sqrt{N}} + g(\varepsilon,N,d) \Big] {\tau_1 + \sigma^2 t \over 2 \theta_d} \Big) \leq 4\varepsilon.
\end{align}
\end{lem}

\begin{proof}[Proof of Lemma~\ref{lem:decoupling_error_bound}]
Recall definitions of the $N$-particle system and the associated approximating $N$ independent OU processes: for $i \in [N]$,
\[
\rd X^{N,i}_{t} = \Theta (\overline{X}^{N}_{t} - X^{N,i}_{t}) \rd t + \sigma \rd W^{i}_{t} \quad \mbox{and} \quad \rd Y^{i}_{t} = -\Theta Y^{i}_{t} \rd t + \sigma \rd W^{i}_{t},
\]
where $(X^{N,i}_{t})$ and $(Y^{i}_{t})$ are driven by the same Brownian motion $(W^{i}_{t})$ and $Y^{i}_0 = X^{N,i}_0$. Averaging the processes $(X^{N,i}_{t})$ over $i \in [N]$, we get $\rd \overline{X}^{N}_{t} = \sigma \rd \overline{W}^{N}_{t} := \sigma N^{-1} \sum_{i=1}^{N} \rd W^{i}_{t},$
i.e., the mean process $(\sqrt{N} \overline{X}^{N}_{t})$ of the $N$ particles has the same law as a standard Brownian motion rescaled by $\sigma$. Combining the last three expressions, we obtain that
\[
\rd (\overline{X}^{N}_{t} - X^{N,i}_{t}+Y^{i}_{t}) = \sigma \rd \overline{W}^{N}_{t} - \Theta (\overline{X}^{N}_{t} - X^{N,i}_{t}+Y^{i}_{t}),
\]
which implies that the difference process $\Delta^{N,i}_{t} = \overline{X}^{N}_{t} - X^{N,i}_{t}+Y^{i}_{t}$ between the $N$-particle system and the decoupled OU processes is an OU process with respect to $(\overline{W}^{N}_{t})$, i.e., we have for all $i \in [N]$,
\begin{equation}
\label{eqn:difference_proc}
\rd \Delta^{N,i}_{t} = -\Theta \Delta^{N,i}_{t} \rd t + \sigma \rd \overline{W}^{N}_{t} \quad \mbox{with} \quad \Delta^{N,i}_0 = \overline{Y}^{N}_0 \sim N(0, N^{-1} D).
\end{equation}
Equation~\eqref{eqn:difference_proc} means that the processes $(\Delta^{N,1}_{t}), \dots, (\Delta^{N,N}_{t})$ are the same mean-reverting OU processes, all driven by $(\overline{W}^{N}_{t})$. Thus $\overline{\Delta}^{N}_{t} = \Delta^{N,i}_{t}$ for all $i \in [N]$, and $N^{-1}\sum_{i=1}^N \int_0^t \rd W^i_s \otimes \Delta^{N,i}_s = \int_0^t \rd \overline{W}_s \otimes \overline{\Delta}^{N}_t$. Now applying~\eqref{eqn:concentration_bound} in Lemma~\ref{lem:concentration_moment_bounds_OU} to the averaged difference process $(\overline{\Delta}^{N}_{t})$ with diffusion parameter $N^{-1/2} \sigma$ and initialization variance $N^{-1} D$, we obtain~\eqref{lem:decoupling_error_bound_1}.

Next we prove~\eqref{lem:decoupling_error_bound_2}. By the triangle inequality and the Cauchy-Schwarz inequality,
\begin{align*}
& \Big\|{1 \over N} \sum_{i=1}^{N} \int_{0}^{t} [ (\overline{X}^{N}_{s} - X^{N,i}_{s}) (\overline{X}^{N}_{s} - X^{N,i}_{s})^T - Y^{i}_{s} {Y^{i}_{s}}^T ] \rd s \Big\| \\
\leq & \Big\|{1 \over N} \sum_{i=1}^{N} \int_{0}^{t} \Delta^{N,i}_s {\Delta^{N,i}_s}^T \rd s \Big\| + \Big\|{1 \over N} \sum_{i=1}^{N} \int_{0}^{t} [ \Delta^{N,i}_s {Y^i_s}^T + Y^i_s {\Delta^{N,i}_s}^T ] \rd s \Big\| \\
\leq & \Big\| \int_{0}^{t} \overline{\Delta}^{N}_s {\overline{\Delta}^{N}_s}^T \rd s \Big\| + 2 \Big( \int_{0}^{t} \| \overline{\Delta}^{N}_s\|^2 \rd s \Big)^{1\over2} \Big( \int_{0}^{t} {1 \over N} \sum_{i=1}^{N}  \| Y^i_s\|^2 \rd s \Big)^{1\over2}.
\end{align*}
Since $(\overline{\Delta}^{N}_t)$ is an OU process, we have $\| \int_{0}^{t} \E [ \overline{\Delta}^{N}_s {\overline{\Delta}^{N}_s}^T ] \rd s \| \leq N^{-1} \int_0^t \|\Sigma_s\| \rd s \leq (\tau_1+\sigma^2 t) (2N\theta_d)^{-1}$, and by~\eqref{eqn:moment_bound},
\[
\Prob \Big( \Big\| \int_{0}^{t} [ \overline{\Delta}^{N}_s {\overline{\Delta}^{N}_s}^T - \E \overline{\Delta}^{N}_s {\overline{\Delta}^{N}_s}^T ] \rd s \Big\| \geq C d {\tau_1 + \sigma^2 t \over N \theta_d} ( \log(d/\varepsilon) + \sqrt{\log(d/\varepsilon)} ) \Big) \leq 2\varepsilon.
\]
So we have with probability at least $1-2\varepsilon$,
\[
\Big\| \int_{0}^{t} \overline{\Delta}^{N}_s {\overline{\Delta}^{N}_s}^T \rd s \Big\| \leq C d g(\varepsilon,1,d)^2  {\tau_1 + \sigma^2 t \over N \theta_d}.
\]
Next, applying~\eqref{eqn:moment_bound_1d} to the process $(\overline{\Delta}^{N}_t)$ and the union bound, we have
\[
\Prob \Big( \max_{j \in [d]} \int_0^t (|\overline{\Delta}^{N}_{sj}|^2 - \E |\overline{\Delta}^{N}_{sj}|^2 ) \rd s \geq C  {\tau_1 + \sigma^2 t \over N \theta_d} ( \log(d/\varepsilon) + \sqrt{\log(d/\varepsilon)} ) \Big) \leq \varepsilon.
\]
This gives
\[
\int_{0}^{t} \| \overline{\Delta}^{N}_s\|^2 \rd s \leq d \, \max_{j \in [d]} \int_0^t |\overline{\Delta}^{N}_{sj}|^2 \rd s \leq C d g(\varepsilon,1,d)^2  {\tau_1 + \sigma^2 t \over N \theta_d},
\]
where the last inequality holds with probability at least $1-\varepsilon$. Similar argument yields
\[
\Prob \Big( {1 \over N} \sum_{i=1}^{N} \int_{0}^{t}  \| Y^i_s\|^2 \rd s \geq C d g(\varepsilon,N,d)^2 {\tau_1 + \sigma^2 t \over \theta_d} \Big) \leq \varepsilon.
\]
Then~\eqref{lem:decoupling_error_bound_2} follows from putting all pieces together.
\end{proof}

%%%%%%%%%%%%%%%%%%%%%%%%%%%%%%%%%%%%%%%%%%%%%%%%%%%%%%%%%%%%%%%%%%%
%%                                                               %%
%% Use the two commands below for producing your bibliography    %%
%% with bibtex, then comment again the commands and include the  %%
%% content of the .bbl file in this file below the commands.     %%
%%                                                               %%
%%%%%%%%%%%%%%%%%%%%%%%%%%%%%%%%%%%%%%%%%%%%%%%%%%%%%%%%%%%%%%%%%%%

\bibliographystyle{amsplain}
\bibliography{interacting_particle_systems}

\providecommand{\bysame}{\leavevmode\hbox to3em{\hrulefill}\thinspace}
\providecommand{\MR}{\relax\ifhmode\unskip\space\fi MR }
\providecommand{\MRhref}[2]{%
  \href{http://www.ams.org/mathscinet-getitem?mr=#1}{#2}
}
\providecommand{\href}[2]{#2}
\begin{thebibliography}{10}

\bibitem{Bongini2017}
Mattia Bongini, Massimo Fornasier, Markus Hansen, and Mauro Maggioni,
  \emph{Inferring interaction rules from observations of evolutive systems i:
  The variational approach}, Mathematical Models and Methods in Applied
  Sciences \textbf{27} (2017), no.~5, 909--951.

\bibitem{Brunton3932}
Steven~L. Brunton, Joshua~L. Proctor, and J.~Nathan Kutz, \emph{Discovering
  governing equations from data by sparse identification of nonlinear dynamical
  systems}, Proceedings of the National Academy of Sciences \textbf{113}
  (2016), no.~15, 3932--3937.

\bibitem{MR2807761}
Peter B{\"u}hlmann and Sara van~de Geer, \emph{Statistics for high-dimensional
  data}, Springer Series in Statistics, Springer, Heidelberg, 2011, Methods,
  theory and applications. \MR{2807761 (2012e:62006)}

\bibitem{ChenXuWu2013}
Xiaohui Chen, Mengyu Xu, and Wei~Biao Wu, \emph{{Covariance and precision
  matrix estimation for high-dimensional time series}}, The Annals of
  Statistics \textbf{41} (2013), no.~6, 2994 -- 3021.

\bibitem{PhysRevLett.96.104302}
M.~R. D'Orsogna, Y.~L. Chuang, A.~L. Bertozzi, and L.~S. Chayes,
  \emph{Self-propelled particles with soft-core interactions: Patterns,
  stability, and collapse}, Phys. Rev. Lett. \textbf{96} (2006), 104302.

\bibitem{EyalYangBahar2006}
Eran Eyal, Lee-Wei Yang, and Ivet Bahar, \emph{{Anisotropic network model:
  systematic evaluation and a new web interface}}, Bioinformatics \textbf{22}
  (2006), no.~21, 2619--2627.

\bibitem{HalilogluBaharErman1997}
Turkan {Haliloglu}, Ivet {Bahar}, and Burak {Erman}, \emph{{Gaussian Dynamics
  of Folded Proteins}}, Phys. Rev. Lett. \textbf{79} (1997), no.~16,
  3090--3093.

\bibitem{HuangLiuLu2019}
Hui Huang, Jian-Guo Liu, and Jianfeng Lu, \emph{Learning interacting particle
  systems: Diffusion parameter estimation for aggregation equations},
  Mathematical Models and Methods in Applied Sciences \textbf{29} (2019),
  no.~1, 1--29.

\bibitem{Kutoyants2004}
Yury~A. Kutoyants, \emph{Statistical inference for ergodie diffusion
  processes}, Springer-Verlag London Ltd., 2004.

\bibitem{LiLuMaggioniTangZhang2019}
Zhongyang Li, Fei Lu, Mauro Maggioni, Sui Tang, and Cheng Zhang, \emph{On the
  identifiability of interaction functions in systems of interacting
  particles}, arXiv:1912.11965 (2019).

\bibitem{Lu14424}
Fei Lu, Ming Zhong, Sui Tang, and Mauro Maggioni, \emph{Nonparametric inference
  of interaction laws in systems of agents from trajectory data}, Proceedings
  of the National Academy of Sciences \textbf{116} (2019), no.~29,
  14424--14433.

\bibitem{McKean1966}
H.~P. McKean, \emph{A class of markov processes associated with nonlinear
  parabolic equations}, Proceedings of the National Academy of Sciences
  \textbf{56} (1966), no.~6, 1907--1911.

\bibitem{MeiE7665}
Song Mei, Andrea Montanari, and Phan-Minh Nguyen, \emph{A mean field view of
  the landscape of two-layer neural networks}, Proceedings of the National
  Academy of Sciences \textbf{115} (2018), no.~33, E7665--E7671.

\bibitem{Mogilner:1999aa}
Alexander Mogilner and Leah Edelstein-Keshet, \emph{A non-local model for a
  swarm}, Journal of Mathematical Biology \textbf{38} (1999), no.~6, 534--570.

\bibitem{MotschTadmor2014}
Sebastien Motsch and Eitan Tadmor, \emph{Heterophilious dynamics enhances
  consensus}, SIAM Review \textbf{56} (2014), no.~4, 577--621.

\bibitem{RevuzYor1991}
Daniel Revuz and Marc Yor, \emph{Continuous martingales and brownian motion},
  Springer-Verlag, 1991.

\bibitem{SchaefferTranWard2018}
Hayden Schaeffer, Giang Tran, and Rachel Ward, \emph{Extracting sparse
  high-dimensional dynamics from limited data}, SIAM Journal on Applied
  Mathematics \textbf{78} (2018), no.~6, 3279--3295.

\bibitem{Sznitman1991}
Alain-Sol Sznitman, \emph{Topics in propagation of chaos}, Ecole d'Et{\'e} de
  Probabilit{\'e}s de Saint-Flour XIX --- 1989 (Berlin, Heidelberg) (Paul-Louis
  Hennequin, ed.), Springer Berlin Heidelberg, 1991, pp.~165--251.

\bibitem{Topaz:2006aa}
Chad~M. Topaz, Andrea~L. Bertozzi, and Mark~A. Lewis, \emph{A nonlocal
  continuum model for biological aggregation}, Bulletin of Mathematical Biology
  \textbf{68} (2006), no.~7, 1601.

\bibitem{Vershynin2018_Cambridge}
Roman Vershynin, \emph{{High-Dimensional Probability: An Introduction with
  Applications in Data Science}}, Cambridge Series in Statistical and
  Probabilistic Mathematics, Cambridge University Press, 2018.

\bibitem{WEN2016237}
Jianghui Wen, Xiangjun Wang, Shuhua Mao, and Xinping Xiao, \emph{Maximum
  likelihood estimation of mckean-vlasov stochastic differential equation and
  its application}, Applied Mathematics and Computation \textbf{274} (2016),
  237 -- 246.

\end{thebibliography}

% add below the content of your .bbl file produced by bibtex.

%\begin{thebibliography}{99}
%
%\bibitem{doob} Doob, J. L.: Heuristic approach to the Kolmogorov-Smirnov
%  theorems. \emph{Ann. Math. Statistics} \textbf{20}, (1949), 393--403.
%  \MR{0030732}
%
%\bibitem{gnekol} Gnedenko, B. V. and Kolmogorov, A. N.: Limit distributions for
%  sums of independent random variables. Translated and annotated by K. L.
%  Chung. With an Appendix by J. L. Doob. \emph{Addison-Wesley}, Cambridge,
%  1954. ix+264 pp. \MR{0062975}
%
%\bibitem{ito} It\^o, K.: Multiple Wiener integral. \emph{J. Math. Soc. Japan}
%  \textbf{3}, (1951), 157--169. \MR{0044064}
%
%\bibitem{levy} L\'evy, P.: Sur certains processus stochastiques homog\`enes.
%  \emph{Compositio Math.} \textbf{7}, (1939), 283--339. \MR{0000919}
%
%\bibitem{grisha} Perelman, G.: The entropy formula for the Ricci flow and its
%  geometric applications, \ARXIV{math.DG/0211159}
%
%\bibitem{smisch} Smirnov, S. and Schramm, O.: On the scaling limits of planar
%  percolation, \ARXIV{1101.5820}
%
%\end{thebibliography}

%%%%%%%%%%%%%%%%%%%%%%%%%%%%%%%%%%%%%%%%%%%%%%%%%%%%%%%%%%%%%%%%%%%
%%                                                               %%
%% You may add acknowledgments (optional).                       %%
%%                                                               %%
%%%%%%%%%%%%%%%%%%%%%%%%%%%%%%%%%%%%%%%%%%%%%%%%%%%%%%%%%%%%%%%%%%%

\ACKNO{Research was supported in part by NSF CAREER Award DMS-1752614 and a Simons Fellowship. Part of this research was carried out in the Institute for Data, System, and Society (IDSS) at Massachusetts Institute of Technology. The author would like to thank Philippe Rigollet (MIT) and Yun Yang (UIUC) for helpful comments.}

%%%%%%%%%%%%%%%%%%%%%%%%%%%%%%%%%%%%%%%%%%%%%%%%%%%%%%%%%%%%%%%%%%%
%%                                                               %%
%% You have reached the end of your document.                    %%
%%                                                               %%
%%%%%%%%%%%%%%%%%%%%%%%%%%%%%%%%%%%%%%%%%%%%%%%%%%%%%%%%%%%%%%%%%%%

\end{document}